\documentclass[11pt, oneside]{article}

\def\ifcover{\ifx34}
\def\widebar{\overline}

\usepackage{graphicx}
\usepackage{amssymb}
\usepackage{amsmath}
\usepackage{amsthm}
\usepackage{indentfirst}

\newcommand\norm[1]{\left\lVert#1\right\rVert}
\newcommand\floor[1]{\lfloor#1\rfloor}
\def\IN{{\mathbb N}}

\def\IR{{\mathbb R}}
\def\P{{\bf P}}
\def\E{{\bf E}}
\def\one{{\bf 1}}
\def\X{\mathcal{X}}

\def\Leb{{\rm Leb}}
\def\half{{1 \over 2}}
\def\bx{{\bf x}}
\def\by{{\bf y}}

\begin{document}

\ifcover

\centerline{\Large\bf Convergence Rates of Attractive-Repulsive MCMC Algorithms}

\bigskip \centerline{by}

\medskip \centerline{Yu Hang Jiang, Tong Liu, Zhiya Lou, Jeffrey S. Rosenthal,}
\medskip \centerline{Shanshan Shangguan, Fei Wang, and Zixuan Wu}

\bigskip\bigskip\noindent
Affiliation of all of the authors:

Department of Statistical Sciences, University of Toronto,
700 University Avenue, 9th Floor,
Toronto, Ontario,
Canada M5G 1X6

\bigskip\bigskip\noindent
Corresponding author:

Jeffrey S. Rosenthal,
Email jeff@math.toronto.edu,
Phone 416-978-4594,
Fax 416-978-5133.

\newpage

\fi

\centerline{\Large\bf Convergence Rates of Attractive-Repulsive MCMC Algorithms}

\bigskip \centerline{by (in alphabetical order)}

\medskip \centerline{Yu Hang Jiang, Tong Liu, Zhiya Lou, Jeffrey S. Rosenthal,}
\medskip \centerline{Shanshan Shangguan, Fei Wang, and Zixuan Wu}

\medskip \centerline{\sl Department of Statistical Sciences, University of Toronto}

\medskip\centerline{(December, 2020; last revised \today)}

\bigskip
\begin{quote}
\noindent \bf Abstract: \rm\
We consider MCMC algorithms for certain particle systems
which include both attractive
and repulsive forces, making their convergence analysis challenging.
We prove that a version of these algorithms on a bounded state space
is uniformly ergodic with explicit quantitative convergence rate.
We also prove that a version on an unbounded state space is
still geometrically ergodic, and then use the method of shift-coupling
to obtain an explicit quantitative bound on its convergence rate.
\end{quote}

\bigskip
\section{Introduction}

Markov Chain Monte Carlo (MCMC) algorithms are
an indispensable tool for researchers and scientists across a wide
spectrum of fields, ranging from machine learning and Bayesian
inference to systems biology and mathematical finance,
to sample from complicated distributions in high dimensions.
When running MCMC, one important question is the number of
steps the Markov chain requires to converge. There are various
approaches to analyzing this difficult problem. In this paper, we
describe a challenging MCMC example,
and show ways of deriving a quantitative mathematical
bound using techniques related to coupling.

\subsection{Background about MCMC}

Markov Chain Monte Carlo (MCMC) algorithms such as the
Metropolis-Hastings algorithm \cite{metropolis1953,hastings}
and the Gibbs sampler \cite{geman1984stochastic,gelfand1990sampling}
have become extremely popular in statistics. They provide a
feasible way to sample from complicated probability distributions
in high dimensions, and play a crucial role in Bayesian inference
as posterior distributions are usually too complicated to compute
analytically. Moreover, the application of MCMC algorithms is not limited
to statistical contexts.  Indeed, the Metropolis algorithm, one
of the most popular MCMC algorithms, arose in physics and was designed
to simulate the behavior of large systems of interacting particles
\cite{metropolis1953}. MCMC algorithms were then widely
applied in computational physics~\cite{barker1965monte,
doi:10.1146/annurev-astro-082214-122339}.
They are now an indispensable tool for researchers and
scientists in many other fields, including computer science
\cite{inproceedings,andrieu2003introduction}, systems biology
\cite{valderrama2019mcmc,whidden2015quantifying}, mathematical
finance \cite{korteweg2011markov,jasra2011sequential}, and more
(e.g.~\cite{geyer1992,brooks2011}).

Specifically, suppose we are given a possibly-unnormalized
density function $\pi(\cdot)$ on a state space $\mathcal{X}$,
e.g.\ a posterior
density in Bayesian statistics. Then, the posterior mean of any
functional $f$ is given by $$\pi(f) \ = \  \frac{\int_{\mathcal{X}}f(x)
\pi(x) dx}{ \int_{\mathcal{X}} \pi(x) dx}.$$ In most cases,
it is infeasible to directly compute this integral (either analytically
or numerically), especially when $\mathcal{X}$ is high-dimensional and
$\pi(\cdot)$ is complicated. An alternative way is to repeatedly sample
from $\pi(\cdot)$, and estimate $\pi(f)$ by the sample average.
However, if $\pi(\cdot)$ is complicated, then it may be
impossible even to draw samples directly from $\pi(\cdot)$. MCMC
algorithms were invented to solve this problem. They construct a Markov
chain which can be easily run on a computer,
which has $\pi(\cdot)$ as its stationary distribution.
It follows under mild conditions that
if we run the Markov chain for a long
time, the distribution of $X_n$ will converge to $\pi(\cdot)$.

In this paper, we will focus on the Metropolis-Hastings algorithm,
one of the simplest and most well-known MCMC algorithms.
Let $\pi(\cdot)$ be an unnormalized density function on
$\mathcal{X}$, and let $q(x, \cdot)$ be an unnormalized density for each
$x \in \mathcal{X}$ . The Metropolis-Hastings Algorithm proceeds as
follows. First we choose some $X_0$ from some initial distribution
$\mu(\cdot)$. Then, for $n=0,1,2,\ldots$,
given $X_n$, we generate a proposal $Y_{n+1} \sim q(X_n,
\cdot)$. With probability $\alpha(X_n, Y_{n+1})$ we set $X_{n+1}
= Y_{n+1}$ where $$\alpha(x, y) = \min \left\{1, \frac{\pi(y)q(y,
x)}{\pi(x)q(x, y)}\right\}$$ is the acceptance rate; otherwise we set
$X_{n+1} = X_n$.
This acceptance probability is chosen precisely to make the Markov chain
reversible with respect to $\pi(\cdot)$, from which it follows that
$\pi(\cdot)$ is a stationary distribution, and under mild conditions the
chain will converge in distribution to $\pi(\cdot)$
\cite{metropolis1953,hastings}.

The knowledge that MCMC will eventually converge to $\pi(\cdot)$
raises the question of how long it takes to converge. There are
various approaches to analyzing this problem. One widely-used method
is to apply diagnostic tools to the output produced by the algorithm
\cite{10.2307/2246093,doi:10.1080/10618600.1998.10474787,
10.2307/2291683}. For example, we can monitor the ergodic averages
of selected scalar quantities of interest (e.g.\ first and second
moments). Another popular
approach is to theoretically derive a bound in terms of the total
variation distance \cite{minor1995,10.2307/2337435,10.2307/3448486},
though this usually involves difficult calculations and the resulting
bounds are often quite conservative. In this paper, we describe
a challenging MCMC example, show ways of deriving a quantitative
mathematical bound using techniques related to coupling, and
compare our theoretical results to diagnostic bounds
from actual computer simulations.

\subsection{The Attractive-Repulsive Model}

We shall focus on the following model.
Suppose we have \(n\) particles randomly located in the  $\mathbb{R}^2$
plane (so the state space $\mathcal{X} = \mathbb{R}^{2n}$), and the
unnormalized density of each configuration is given by
\begin{equation}\label{mainpi}
\pi(x) = \exp\Big(- \Big[c_1 \sum_{i = 1}^n||x_i|| + c_2 \sum_{i < j}
||x_i - x_j||^{-1}\Big]\Big),
\end{equation}
where $c_1$, $c_2$ are
positive constants and  $||\cdot||$ is the usual Euclidean ($L^2$)
norm on $\mathbb{R}^2$.  Since the density is fairly complicated, it
is hard to compute expected values with respect to this distribution,
such as the average distance of the particles
to the origin. Therefore, a more feasible solution is to simulate
this distribution using an MCMC algorithm. We shall use
componentwise versions of the Metropolis-Hastings algorithm
\cite{metropolis1953,hastings},
in which the multiple particles are updated one
at a time in a sequential order, each with a proposal followed by an
accept/reject step.
(For a graphical illustration of this algorithm on these
densities, see \cite{pointproc.js}.)
By running the algorithm for many iterations, we can approximately
sample from $\pi$, and thus find good estimates of its expected values.

The density function~\eqref{mainpi} is designed so the first
summation ``pulls" the particles towards the origin, while the second
summation ``pushes" them away from each other. Hence, we call this
an \textit{attractive-repulsive} particle system.  The combination of
attractive and repulsive forces mean that the MCMC algorithm does not
satisfy simple monotonicity or other properties which would simplify its
convergence analysis, so that more careful techniques are required.
Nevertheless, for certain special cases of this density, we will derive
both qualitative and quantitative convergence bounds herein.

We note that there is a long history of using MCMC to study interacting
particle models.  For example,
Alder and Wainwright~\cite{alder} used Monte Carlo to simulate
the dynamics of molecules;
Hammersley~\cite{hamm1972} and Liggett~\cite{ligg1978} applied stochastic
atomic lattice models to solid-state physics particle systems;
Speagle~\cite[Section~8]{speagle2020conceptual} studied
a purely attractive model (where a
particle is more likely to move inwards than outwards) using a Metropolis
algorithm with Gaussian proposal distributions;
and Krauth~\cite{krauth2021eventchain}
used local non-reversible MCMC algorithms to
simulate dynamic hard-spheres.
The model~\eqref{mainpi} is similar in spirit to these other dynamics,
though it was chosen primarily for illustrative purposes (e.g.\
it is not stochastically monotone; see below).

\subsection{Background about Minorization and Drift Conditions}

We are interested in bounding the {\it total variation distance}
$$
\|P^n(x,\cdot) - \pi(\cdot)\|
\ := \ \sup_{S \subseteq \X} |P^n(x,S) - \pi(S)|
\ = \ \sup_{S \subseteq \X} |P(X_n \in S \, | \, X_0=x) - \pi(S)|
$$
between the
\(n\)-step distribution $P^n(x,\cdot)$
and the stationary distribution \(\pi(\cdot)\)
of a Markov chain, where the supremum is taken over all
measurable subsets $S$.  One method involves coupling via
minorization and drift conditions. A Markov chain with a
state space $\mathcal{X}$ and transition probabilities
$P(x, \cdot)$ satisfies a {\it minorization condition} if there is a
measurable subset $C\subseteq \mathcal{X}$, a probability measure $Q$
on $\mathcal{X}$, a constant $\epsilon >0$, and a positive integer
$n_0$, such that
\begin{equation}\label{minorcond}
P^{n_0}(x,\cdot) \geq \, \epsilon \, Q(\cdot), \quad x\in C.
\end{equation}
We call such $C$ a {\it small set}, and refer to it
$(n_0,\epsilon, Q)$-small.
In particular, if $C=\mathcal{X}$ (i.e., $C$
is the entire state space), then we say the Markov chain satisfies a
{\it uniform} minorization condition, also referred to as
{\it Doeblin's condition} (see \cite{doeblin1938expose}).
It then follows (see e.g.~\cite{MeynTweedie,probsurv})
that the chain is {\it uniformly ergodic}, i.e.\
there are fixed $\rho <1$ and $M < \infty$ such that
$$
\norm{P^n(x,\cdot)-\pi(\cdot)} \ \leq \ M \, \rho^n,
\quad n\in\IN,
\quad x\in\X,
$$
and in fact a precise convergence bound is available:

\paragraph{Proposition 1:}
If a Markov chain with stationary distribution $\pi(\cdot)$
has the property that the entire state space
$\mathcal{X}$ is $(n_0,\epsilon, Q)$-small, then the chain
is uniformly ergodic, with
$$
\norm{P^n(x,\cdot) - \pi(\cdot)} \leq
(1-\epsilon)^{\floor{\frac{n}{n_0}}}
\, ,
\quad n\in\IN
\, .
$$
\medskip

In Section~2, we prove uniform ergodicity for a bounded version of our
algorithm.
Unfortunately, many Markov chains are not uniformly ergodic.
A Markov chain with stationary distribution
\(\pi(\cdot)\) is {\it geometrically ergodic} if 
there are fixed $\rho <1$ and $\pi$-a.e.-finite function $M:\X\to[0,\infty]$
such that
$$
\norm{P^n(x,\cdot) - \pi(\cdot) } \ \leq \ M(x) \, \rho^n,
\quad n \in \mathbb{N},
\quad x\in\X,
$$
i.e.\ if the multiplier $M$ can depend on the initial state~$x$.
Also, a Markov chain with a small set $C$ satisfies a {\it univariate drift
condition} if there are constants $0<\lambda<1$ and $b<\infty$, and a
$\pi$-a.e.-finite
function $V: \mathcal{X}\rightarrow [1,\infty]$ such that
\begin{equation}\label{driftcond}
PV(x)
\ := \ E[V(X_1) \, | \, X_0=x]
\ \leq \ \lambda \, V(x) + b \, \mathbf{1}_C(x), \quad x \in \mathcal{X}.
\end{equation}
The minorization condition~\eqref{minorcond}
and drift condition~\eqref{driftcond} together
guarantee that the chain is geometrically ergodic
(e.g.\ \cite[Theorem~15.0.1]{MeynTweedie}):

\paragraph{Proposition 2:} If a $\phi$-irreducible, aperiodic
Markov chain with stationary distribution $\pi(\cdot)$
and small set \(C \subset \mathcal{X}\) satisfies
the minorization condition~\eqref{minorcond} for some $n_0\in\IN$
and $\epsilon>0$ and $C \subseteq \X$ and probability measure $Q(\cdot)$ on
$\X$, and the drift condition~\eqref{driftcond} for some
$\pi$-a.e.-finite function $V:\X\to[0,\infty]$ and $\lambda<1$
and $b<\infty$,
then it is geometrically ergodic.
\bigskip\medskip

Geometric ergodicity is a helpful property, since it implies the chain
converges geometrically quickly, and also implies certain other results
such as central limit theorems (see e.g.\ \cite{probsurv}).
We establish it for an unbounded version of our algorithm in Section~3.
Unfortunately, qualitative bounds such as
uniform or geometric ergodicity can still
be quite weak in many cases, and do not necessarily imply that the Markov chain
converges in a short time.
For example, if $\mathcal{X} = \{0, 1\}$, with $X_0=x=1$ and
$$P = \begin{pmatrix} 1& 0  \\ 1-z & z\end{pmatrix}$$
for some fixed $z\in (0, 1)$, then $\pi=(1,0)$, and the chain satisfies
a uniform minorization condition with $\epsilon = 1 - z$ and $Q = (1, 0)$.
So, it is both uniformly and geometrically ergodic, and in fact
$\norm{P^n(x,\cdot) - \pi(\cdot)} \, = \, z^n$.
However, it converges
arbitrarily slowly for $z$ near~1, indicating that geometric ergodicity
does not really imply fast convergence.
Due to these limitations, it is best
to find a {\it quantitative} bound, i.e.\ explicit bounds on
$\norm{P^n(x,\cdot) - \pi(\cdot)}$ which provide a value of $n$ that
guarantees that this distance will be sufficiently small.
We consider this problem for an unbounded version of
our attractive-repulsive processes in Section~4 below.

\subsection{Organisation of the Paper}

This paper is organised as follows.  In Section~2, we consider a version
of our algorithm within a bounded domain, and show that it is uniformly
ergodic by means of an explicit uniform minorization condition. In
Section~3, we expand the state space to all of \(\mathbb{R}^2\), and show
that a version of our algorithm is still geometrically ergodic since
it satisfies an explicit univariate drift condition.  In Section~4,
we discuss the challenges of computing a quantitative convergence
bound for our algorithm, and use a {\it shift coupling} construction
to overcome these problems and obtain an explicit quantitative bound.
In Section~5, we compare our theoretical results to observed convergence
behaviour from actual computer simulations.
In Section~6, we provide proofs of all of the theorems in this paper.

\section{Particles in a Square: Uniform Ergodicity}

In this section, we study the attractive-repulsive particle system
density~\eqref{mainpi} in a compact setting.  Suppose we have $n=3$
particles randomly located in the square $U = [0, 1]^2 \subset
\mathbb{R}^2$, with the particle positions denoted by $\bx = (x_i)_{i
= 1, 2, 3} = (x_{i1}, x_{i2})_{i = 1, 2, 3}$, so the state space
$\mathcal{X} = [0, 1]^6$.

We use a componentwise Metropolis
algorithm with systematic scan, in which we repeatedly
update the $n=3$ particles in order (see e.g.\ \cite{metropolis1953,
brooks2011, pointproc.js}).  Specifically, given a configuration $X_n
= \bx$, we first ``propose'' a new location for the first particle
$x_1$ from the uniform (Lebesgue) measure on $\mathcal{X}$, to obtain a
new particle location $y_1$, and hence a
new proposed
configuration $\by=(y_1,x_2,x_3)$. Then with probability $\alpha(\bx,\by)=
\min\big[1,\frac{\pi(\by)}{\pi(\bx)}\big]$, we ``accept'' this proposal
and update $x_1$ to $y_1$.
Otherwise, we ``reject'' this proposal and leave
the original $x_1$ unchanged. We then similarly update $x_2$ and then $x_3$.
That entire procedure represents one iteration of our algorithm, which we
then repeat $n$ times to obtain a final configuration $X_n$.

For this algorithm, we show (all theorems are proved in
Section~\ref{sec-proofs}):

\medskip\noindent
\textbf{Theorem 1.} The above Markov Chain
(a componentwise Metropolis algorithm with uniform proposals and
systematic scan, for
the unnormalised density~\eqref{mainpi} on $[0,1]^6$ with
$n=3$ particles for some constants $c_1,c_2>0$)
is uniformly ergodic,
and satisfies a uniform minorization condition with $n_0 = 1$
and $\epsilon = (0.48)e^{- c_1 (8.49) - c_2 (19.76)}$.

\bigskip

For example, if $c_1 = c_2 = 1/10$, then we can take $\epsilon = 0.028$. By
Proposition 1, we have
$\|P^n(\bx, \cdot) - \pi(\cdot)\| \le (0.972)^n$.
This proves that after 163 steps, the total variation distance
between the \(n\)-step distribution and the stationary distribution
\(\pi(\cdot)\) of this Markov chain will be within 0.01.

\medskip\noindent\bf Remark. \rm
The above model and algorithm could also be considered for $n>3$ particles,
and the convergence rate could probably be bounded in that case too
by similar methods, but the computations become messier,
so here we stick to $n=3$ particles for ease of analysis.

\section{One particle in $\mathbb{R}^2$: Geometric Ergodicity}

We now extend our state space to the entire $\mathbb{R}^2$ plane,
but with just $n=1$ particle. Specifically, suppose we have a particle
randomly located at $ \mathbb{R}^2$, denoted by $x = (x_1, x_2)$,
with unnormalized density given by
\begin{equation}\label{onepi}
\pi(x) \ = \ e^{- H(x)}, \quad \hbox{\rm where} \quad
H(x) \ = \ ||x|| + \frac{1}{||x||}
\ := \ r_x + \frac{1}{r_x},
\end{equation}
where $r_x := ||x||$ is again the $L^2$ norm.
Note that this model~\eqref{onepi} can be considered to be a special case
of our main model~\eqref{mainpi},
in which $c_1=c_2=1$ and $n=2$, where one particle is at~$x_1 := x$,
and a second particle is always fixed to be at the origin $x_2 := 0$.

We use the following Metropolis-Hastings algorithm on this distribution.
For any $x = (x_1, x_2) \in \mathbb{R}^2$, let
$$
B_x \ = \ \{ z\in\IR^2 : |r_x-1| < \|z\| < r_x+1 \}
\, .
$$
Thus, $B_x$ is an annulus of width $2 \min(r_x,1)$, which contains $x$
unless $r_x<0.5$; see Figure~1.
And, $vol(B_x) = \pi(r_x + 1)^2 - \pi|r_x - 1|^2 = 4\pi r_x$.
We then let the proposal density $q(x,
\cdot)$ be the uniform distribution on $B_x$, i.e.\
$$
q(x, dy) \ = \ \one_{B_x}(y) \  \frac{dy}{4\pi r_x}
\, , \quad x,y \in R^2
\, .
$$

\begin{figure}[ht]
\centering
\includegraphics[scale=0.35]{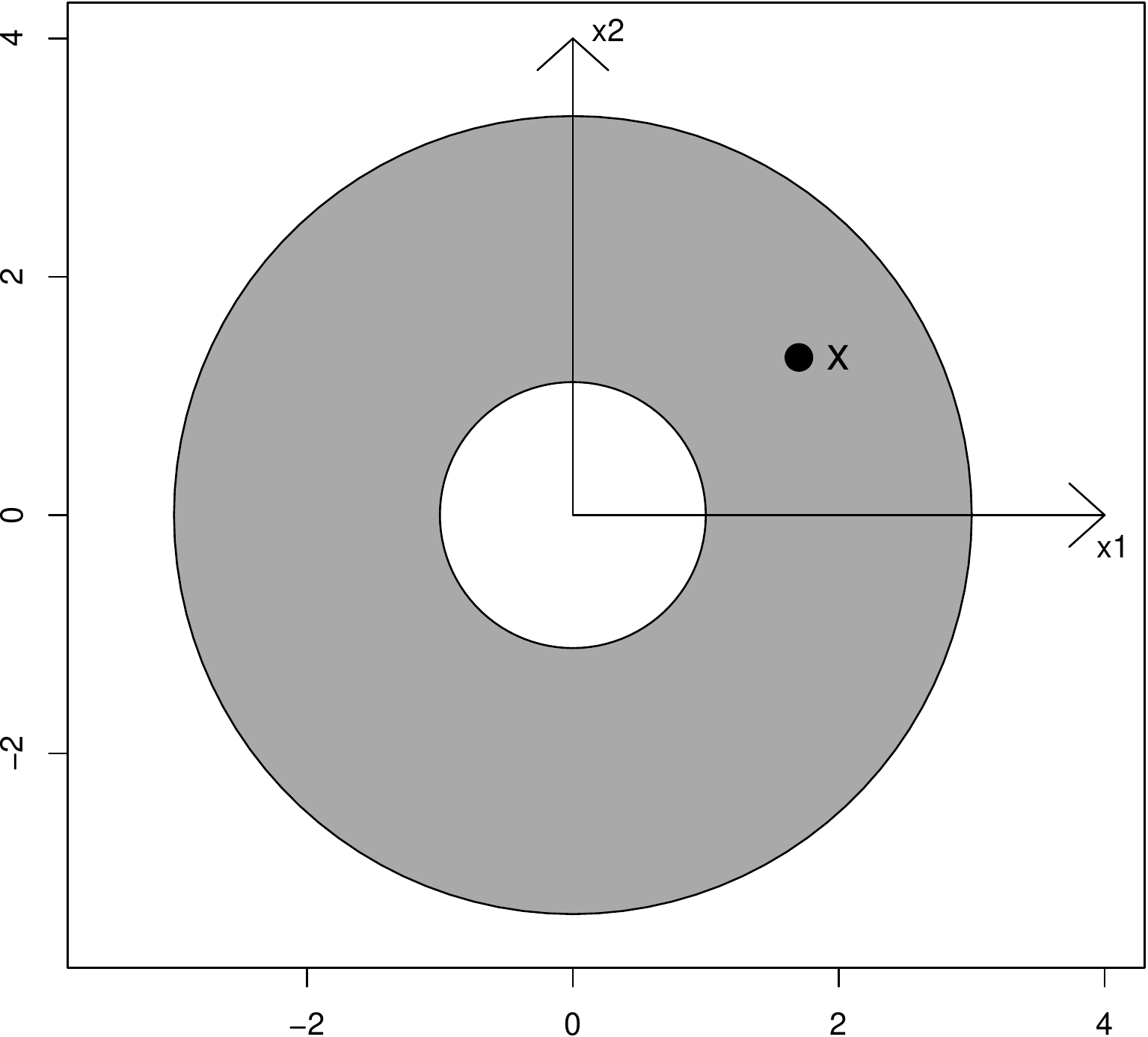}
\ \ \includegraphics[scale=0.35]{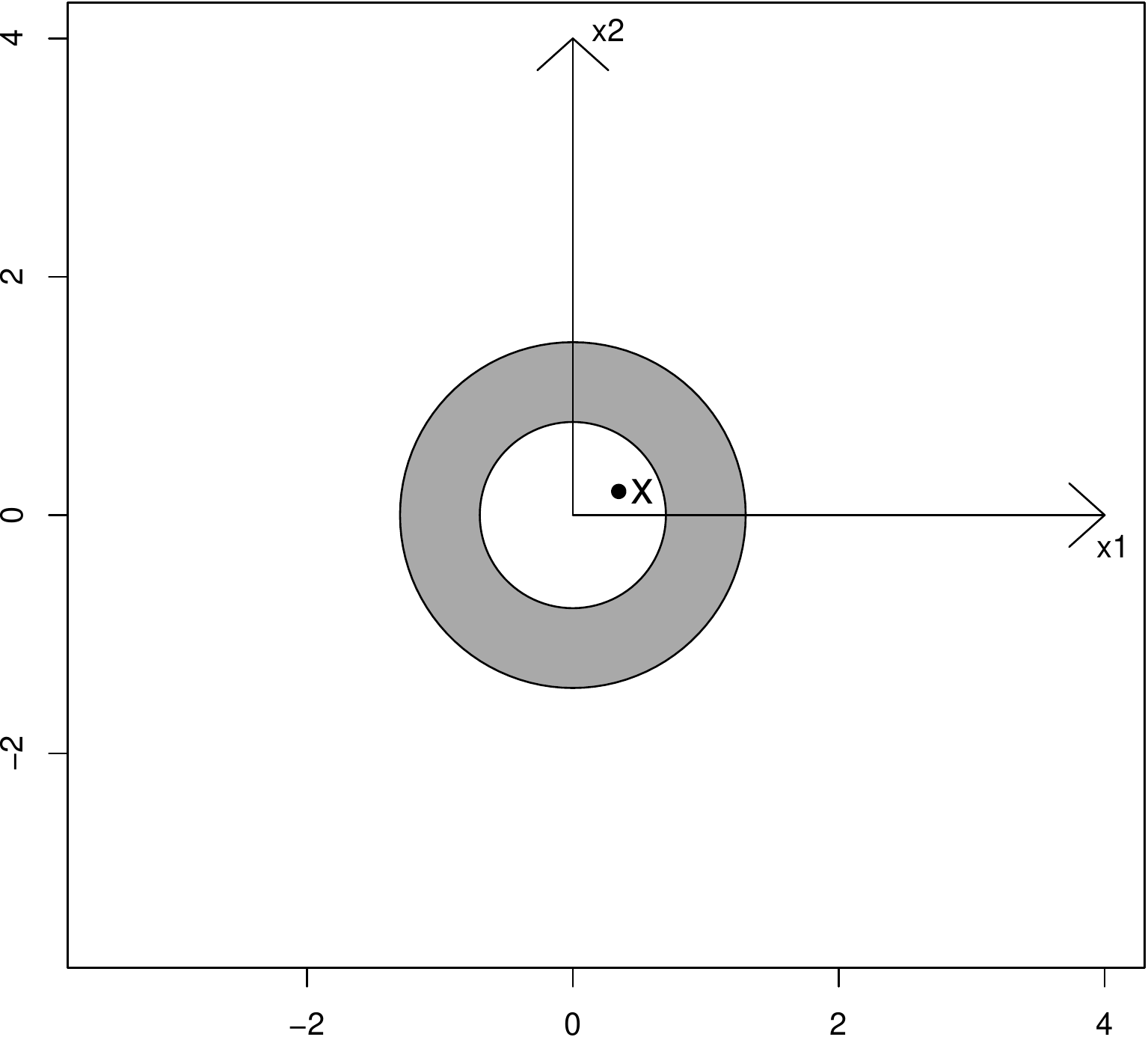}
\caption{Illustration of the region $B_x$ for the model of Section~3,
in two different cases: when $r_x=2$ (left) or $r_x=0.3$ (right).}
\centering 
\end{figure}

Note that $y \in B_x$ if and only if $x \in B_y$ (since for $r_x,r_y
\le 1$ this is equivalent to $r_x+r_y<1$; and for $r_x<1<r_y$
or $r_y<1<r_x$ this is
equivalent to $\min[r_x,r_y] < \max[r_x,r_y]+1$;
and for $r_x,r_y \ge 1$ this is equivalent
to $|r_x-r_y|<1$).  Hence, these $q(x,dy)$ are valid proposal distributions
for a Metropolis-Hastings algorithm.  The corresponding acceptance rate is
$$
\alpha(x, y)
\ = \ \min\bigg\{1, \ \frac{\pi_u(y) \, q(y, x)}{\pi_u(x) \, q(x, y)} \bigg\}
\ = \ \min\bigg\{1, \ \frac{e^{H(x)} \, r_x}{e^{H(y)} \, r_y}\bigg\}.
$$
This algorithm is somewhat related to the algorithm of Section~2, except
with just one particle to move so there is no ``scan''
of different particles, and
with a more complicated proposal distribution since the state space $\X$
is unbounded.

For the above algorithm, we shall prove the following quantitative conditions:

\medskip\noindent
\textbf{Theorem 2.} The Markov chain constructed above
(a Metropolis-Hastings algorithm 
with proposals uniform on $B_x$, for
the unnormalised density~\eqref{onepi} on $\IR^2$ with one particle)
satisfies:\\
 \textbf{(a)} the minorization condition
$$P^2(x, \cdot) \ \geq \ (3.5 \times 10^{-5}) \, Q(\cdot),\quad x \in C,$$
\indent for some $Q(\cdot)$,
where $C = \{x \in \mathbb{R}^2, \ \frac{1}{4} \le ||x|| \le 4\}
\subseteq \mathcal{X}$.\\
 \textbf{(b)} the univariate drift condition
$$PV(x) \ \le \ 0.995 \, V(x) + (e^{2.7}-0.995) \, \mathbf{1}_C,
\quad x \in \mathcal{X},$$ 
\indent where $V(x) = e^{\frac{1}{2}H(x)}$.
Furthermore, $\sup_{x\in C} PV(x) \le e^{2.7}$.
\hfil\break
In particular, by Proposition~2, this chain is geometrically ergodic.

\medskip\noindent\bf Remark. \rm
The above model and algorithm could also be considered for $n>1$ particles,
and it is possible that
the convergence rate could probably be bounded in that case too,
but the analysis becomes much more challenging,
so here we stick to just~1 particle for ease of analysis.

\section{Quantitative Bounds and Shift Coupling}

We next consider quantitative bounds for the algorithm in the previous
section.  There are many potential ways to obtain quantitative bounds
for MCMC algorithms. However, not all methods are feasible for our
attractive-repulsive process.

One common approach uses minorization conditions and bivariate drift conditions
(e.g.\ \cite{minor1995,honest}).
Theorem~3 already provides a minorization condition
and a univariate drift condition, and there are ways to derive a bivariate
drift condition from a univariate one if certain conditions are
satisfied (see e.g.\ Proposition~11 of \cite{probsurv}). However,
to obtain a bivariate drift condition for our processes, we would
have to prove a multi-step minorization condition on a much larger subset,
which would be very challenging and lead to extremely weak bounds.

\def\stochle{\preccurlyeq}

Alternatively, minorization and univariate drift conditions give good
quantitative bounds for Markov chains which are {\it stochastically
monotone},
meaning that there is some stochastic ordering $\stochle$ on $\X$ which is
probabilistically preserved
\cite{Daley1968,cohn1983,lund1996,roberts_tweedie_2000}).
More formally, $P(x_1, B_y) \ge P(x_2,B_y)$ for all $x_1,x_2,y\in\X$
with $x_1 \stochle x_2$, where $B_y = \{ z\in\X : z \stochle y \}$.
Indeed, if we considered a purely attractive version of our model,
by setting $c_2=0$ in~\eqref{mainpi}, then our Markov chain would indeed
be stochastic monotone under the partial order defined
by $x \stochle y$ if and only if $\|x\| \le \|y\|$.  However,
with $c_1,c_2>0$,
the attractive-repulsive nature of our model~\eqref{mainpi}
seems to preclude any stochastic monotonicity condition, so
the improved convergence bounds for stochastically monotone Markov chains
cannot be applied.

Instead, we shall use a particular coupling
method called {\it shift coupling} \cite{aldous1993,shift} to derive
a quantitative bound for the particle system. This
construction only requires a univariate drift condition (not a bivariate one),
and does not require aperiodicity.
In the shift coupling construction, just like ordinary coupling, we
will jointly define two Markov chains to obtain a bound on the rate
of convergence. The key point in which shift coupling differs from the
ordinary method is that we allow the chains to couple at different times.

Let $P(\cdot, \cdot)$ be the transition probabilities for a Markov chain
on a state space $\mathcal{X}$.  Assume the chain is
$\phi$-irreducible, with stationary distribution $\pi(\cdot)$.
Let $\{X_k\}_{k = 0}^{\infty}$ and $\{X_k'\}_{k = 0}^{\infty}$
be two different copies of the chain, defined jointly. Suppose
$T$ and $T'$ are two random variables taking values in $\mathbb{Z}_{\geq
0} \cup \{\infty\}$, such that for any non-negative integer $n$, $X_{T +
n} = X_{T' + n}$.
Ordinary coupling requires $T = T'$,
but shift coupling allows the two Markov chains to become equal at
different times, thus making it easier for the chains to couple.
We can then combine this shift-coupling bound with minorization and
univariate drift conditions, leading to the following
(which generalizes Theorem~4 of \cite{shift} to the case $n_0>1$):

\medskip
\noindent \textbf{Theorem 3:} Suppose a Markov chain
on a state space $\mathcal{X}$, with initial distribution $\nu(\cdot)$,
transition probabilities $P(\cdot, \cdot)$,
and stationary distribution $\pi(\cdot)$, satisfies
the minorization condition~\eqref{minorcond} for some $n_0\in\IN$
and $\epsilon>0$ and $C \subseteq \X$ and probability measure $Q(\cdot)$ on
$\X$, and the drift condition~\eqref{driftcond} for some
$\pi$-a.e.-finite function $V:\X\to[0,\infty]$ and $\lambda<1$
and $b<\infty$,
such that $C= \{x\in \mathcal{X} : V(x) \leq d\}$ for some fixed $d\geq 0$.
Then setting $A := \sup_{x\in C}E(V(X_1)|X_0=x)$
(so $A \leq \lambda d+b$),
for any $0<r<1$ such that $\lambda^{(1-n_0r)} A^r < 1$, we have
\begin{align*}
&\bigg\|\frac{1}{n}\sum_{k=1}^n P(X_k\in \cdot)-\pi(\cdot)\bigg\| \cr
    &\leq \frac{1}{n} \left[\frac{2(1-\epsilon)^r}{1-(1-\epsilon)^r}+
    \frac{\lambda^{-n_0+1-n_0 r} A^r}{1-\lambda^{1-n_0 r} A^r}\left
    (E_\nu(V)+\frac{b}{1-\lambda}\right)\right].
\end{align*}

\bigskip

We now apply this shift-coupling bound to the
attractive-repulsive particle systems of Section~3.
By Theorem~2, we can take
$\epsilon = 3.5 \times 10^{-5}$, $n_0 = 2$,
$\lambda = 0.995$,
$b = e^{2.7} - 0.995$,
$d = e^{17/8}$,
and $A=e^{2.7}$.
Assume the chain starts from the point $(1,0)$,
so $E_\nu(V) = V((1,0)) = e^{\half(1+{1 \over 1})} = e$.
Choosing $r = 0.0016$,
we have $\lambda^{(1 - n_0 r)}A^r \doteq 0.9993 < 1$, and
we compute from Theorem~3 that
$$
\bigg\|\frac{1}{n}\sum_{k=1}^n \P(X_k\in.)-\pi(.)\bigg\| \ \leq \
\frac{39,900,000}{n}.
$$
This bound is certainly far from tight.  However, it does show that
shift-coupling can provide explicit quantitative bounds
on the distance to stationarity, even for the attractive-repulsive
processes that we consider herein.

Finally, we note that the left-hand side of the bound in Theorem~3
differs from the conventional total variation distance between
the \(n\)-step distribution and the stationary distribution. This raises
the question of the meaning of the quantity we are bounding.
An interpretation is given by the following result.

\medskip\noindent
\textbf{Theorem~4:}
Let $\{X_k\}$ be a Markov chain on a state space $\X$,
with transition probabilities $P(\cdot,\cdot)$ and stationary
distribution $\pi(\cdot)$.
For $n \in \mathbb{N}$ and measurable $S \subseteq \X$, let $F_n(S)
:= \E\big[ {1 \over n} \ \# \{ i : 1 \le i \le n , \ X_i \in S\} \big]$
be the expected fraction of time from $1$ to $n$
that the chain is inside $S$.  Then 
$$
\sup_{S\subseteq\X} |F_n(S) - \pi(S)|
= \bigg\|\frac{1}{n}\sum_{k=1}^n \P(X_k\in \cdot)-\pi(\cdot)\bigg\|
\le \frac{1}{n}\sum_{k=1}^n \bigg\|\P(X_k \in \cdot) - \pi(\cdot)\bigg\|.
$$

\medskip
Theorem~4 provides context for Theorem~3.  It shows that the bound of
Theorem~3 in turn provides an upper bound on the difference between the
expected occupation fraction of $S$ and the target probability $\pi(S)$,
uniformly over choice of subset $S$.  So, if the bound is small, then the
chain spends approximately the target fraction of time in every subset, on
average.

Theorem~4 also gives us a way to relate the shift coupling result
to more conventional results.  In particular, note that $||\P(X_k \in
\cdot) - \pi(\cdot)||$ is the usual total variation distance
discussed in previous sections.  Hence, $\frac{1}{n}\sum_{k=1}^n
||\P(X_k \in \cdot) - \pi(\cdot)||$ is the average of the total
variation distances between the $k$-step distribution and the stationary
distribution, averaged over $k=1,2,\ldots,n$.
However, due to the inequality,
$|\frac{1}{n}\sum_{k=1}^n \P(X_k\in \cdot)-\pi(\cdot)|$
does {\it not} provide an upper bound for
$\frac{1}{n}\sum_{k=1}^n |\P(X_k \in \cdot) - \pi(\cdot)|$.

\section{Simulations -- Convergence Diagnostics}

In this section, we run the MCMC algorithms discussed in Sections~2
and~3 above, and apply the MCMC convergence diagnostic tools of
\cite{10.2307/2246093,doi:10.1080/10618600.1998.10474787} to estimate
their convergence times by comparing between- and
within-chain variances of multiple runs of the algorithm when starting
from an over-dispersed starting distribution.
We then compare these estimated times with the
theoretical bounds derived in the previous sections.

\subsection{Three particles in a square}

We begin with the model of Section~2, i.e.\
the componentwise Metropolis algorithm with uniform proposals and
systematic scan for
the unnormalised density~\eqref{mainpi} on $[0,1]^6$ with
$n=3$ particles for some constants $c_1,c_2>0$.
We use the uniform distribution on $[0,1]^6$
as our over-dispersed starting distribution.
To proceed, following
\cite{10.2307/2246093,doi:10.1080/10618600.1998.10474787},
we draw $m=5$ initial samples from this
starting distribution, and then
run $m=5$ different chains in parallel, each for $n=60$ iterations.

Our goal is to see if the chain has converged after $n_*=30$ iterations,
i.e.\ if iterations $n_*+1$ through $n$
(i.e., 31 through 60) are approximately in stationarity.
To investigate this, for iterations $n_*+1$ through $n$ and
initial test functional
$$
\psi: \mathbb{R}^6 \to \mathbb{R} \quad {\rm by} \quad \ \psi(\bx) =
\sqrt{x_{11}^2 + x_{12}^2}
+ \sqrt{x_{21}^2 + x_{22}^2}
+ \sqrt{x_{31}^2 + x_{32}^2}
\ ,
$$
we calculate the
between-chain variance $B$ and the within-chain variances $W$:
$$
B \ := \ \frac{n-n_*}{m - 1} \ \sum_{j = 1}^m (\widebar{\psi}_j - \widebar{\psi})^2,
$$
$$
W \ := \ \frac{1}{m} \, \sum_{i = 1}^m s_i^2
\ = \ \frac{1}{m(n-n_* - 1)} \ \sum_{j = 1}^m \
\sum_{t= n_*+1}^n (\psi_{jt}- \widebar{\psi}_j)^2,
$$
where $n = 60$, $n_*=30$,
$m = 5$, $\psi_{jt} = \psi(X_{jt})$
is the value of $\psi$ on the $t^{\rm th}$ iteration of
chain $j$, $\widebar{\psi}_j = {1 \over n-n_*} \sum_{t=n_*+1}^n \psi(X_{jt})$
is the sample mean of $\psi$ in chain $j$
over iterations $n_*+1$ through $n$,
and $\widebar{\psi} = {1 \over m} \sum_{j=1}^m \widebar\psi_j
= {1 \over m(n-n_*)} \sum_{j=1}^m \sum_{t=n_*+1}^n \psi(X_{jt})$
is the mean of the $m$ different $\widebar{\psi}_j$ values
(i.e.\ the mean of all $m(n-n_*)$ post-burn-in simulated values).

In our simulations, we obtained the values
$$
B =   0.4899, \ \ W = 0.19.
$$
We then estimate the target
variance by a weighted average of $B$ and $W$: $$\hat{\sigma}^2 :=
\frac{B}{n} + \frac{n - 1}{n}W = 0.200.$$
We also compute the pooled posterior variance
estimate $$\hat{V} := \hat{\sigma}^2 + \frac{B}{mn} = 0.2033.$$
Finally we compute the potential scale reduction factor (PCRF),
as $$R \ := \ \frac{d + 3}{d + 1} \cdot \frac{\hat{V}}{W} \ = \
1.07,$$ where $d$ is the degrees of freedom of the corresponding
t-distribution (so $(d + 3)/(d + 1) \approx 1$).
This produced the value $R = 1.07$.  Since this value is $< 1.2$,
that fact provides some indication
\cite{10.2307/2246093,doi:10.1080/10618600.1998.10474787}
that the chain might have approximately converged after $n_*=30$ iterations
(though this diagnostic does {\it not} directly estimate the
total variation distance; see Section~6).

We also consider some other test functionals.
Let $$\phi_1: \mathbb{R}^6 \to
\mathbb{R} \quad {\rm by} \quad \ \phi_1(\bx) = x_{11} + x_{12} + x_{21} + x_{22} + x_{31}
+ x_{32},$$
and
$$\phi_2: \mathbb{R}^6 \to \mathbb{R} \quad {\rm by} \quad \ \phi_2(\bx)   =
x_{11} \cdot x_{12} + x_{21} \cdot x_{22} + x_{31} \cdot x_{32}.$$
Following the same steps as above,
we compute the corresponding PCRFs: $$R_1 = 1.092, \ \ R_2 =
1.091.$$
These values are again $<1.2$.  Hence,
these test results all provide some indication
that the chain might have approximately converged
after $n_*=30$ iterations.  If so, then this is
somewhat quicker than the theoretical bound
(163 iterations) derived in Section~2, suggesting that our bound is
overly conservative.
However, there is a clear benefit in having definitive, guaranteed
(though conservative) theoretical bounds, rather than relying on
convergence diagnostics which can sometimes be misleading
(cf.\ \cite{matthews,cowles}).

\subsection{One particle in $\mathbb{R}^2$}

We next consider the model of Section~3,
i.e.\ the Metropolis-Hastings algorithm 
with proposals uniform on $B_x$, for
the unnormalised density~\eqref{onepi} on $\IR^2$ with one particle.
For our over-dispersed starting distribution we take
the uniform distribution on $[-10, 10]^2$.  We
draw $m=10$ samples from it, as the starting states for 10 different chains,
each run for $n=600$ iterations.

Our goal is to see if the chain has converged after $n_*=300$ iterations,
i.e.\ if iterations $n_*+1$ through $n$
(i.e., 301 through 600) are approximately in stationarity.
We then run our $m=10$ different chains in parallel, each for $n=600$
iterations, and investigate iterations $n_*+1$ through $n$.
For our test function, we begin with
$$\psi(x): \mathbb{R}^2 \to \mathbb{R} \quad {\rm by} \quad \ \psi(x) = ||x||.$$
For this function,
we calculate the between-chain variance $B$ and the within-chain
variance $W$ as above, to obtain: $$B= 60.5, \ \ W = 2.048.$$
We then compute
the corresponding pooled variance and PCRF values to be:
$$\hat{V} = 2.283, \ \ R = 1.115.$$
We again have $R<1.2$, which provides some indication that the chain might
have approximately converged $n_*=300$ iterations.

To investigate further, we consider the two additional test functions
$$\phi_1: \mathbb{R}^2 \to \mathbb{R}
\quad {\rm by} \quad \ \phi_1(x) = |x_1| + |x_2|,$$
and
$$\phi_2:
\mathbb{R}^2 \to \mathbb{R} \quad {\rm by} \quad \ \phi_2(x)
= \begin{cases} 1, & 0.5
\le ||x|| < 1.5 \\ 0, & \text{otherwise} \end{cases}$$
For these test functions, we compute
the corresponding PCTF values to be
$$R_1 = 1.115, \ {\rm and} \ R_2 = 1.061.$$
These values are all $<1.2$, so all of these test results
again provide some indication that the chain might have
approximately converged after $n_*=300$ iterations.
Once again, this is much quicker than the overly-conservative theoretical
bounds derived in Section~4 above.
However, there is again potential benefit in having guaranteed
theoretical bounds, rather than just suggestive convergence diagnostics.

\section{Simulations -- Total Variation Distance}

For more direct comparison with our theoretical results, we
now attempt to estimate
the actual total variation distance between the stationary
distribution and the simulated Markov chain distribution after
different numbers of iterations.
Recall~\cite[Proposition~3(b)]{probsurv}
that one of the many equivalent definitions between two
probability distributions $\nu_1(\cdot)$, $\nu_2(\cdot)$ is
$$||\nu_1(\cdot) - \nu_2(\cdot)||_{TV} = \frac{1}{b-a}\sup_{f
: \mathcal{X} \to [a, b]} \left|\int f d \nu_1 - \int f
d\nu_2\right|,$$ where $a < b$ are real numbers. We shall apply
this definition with different choices of functional $f$
to estimate the total variation distance to stationarity.

\subsection{Three particles in a square}

We first consider the model of Section~2, i.e.\
the componentwise Metropolis algorithm with uniform proposals and
systematic scan for
the unnormalised density~\eqref{mainpi} on $[0,1]^6$ with
$n=3$ particles for some constants $c_1,c_2>0$.
We apply different functionals to estimate the total variation
distance. We begin with the functional
$$f: \left[0, 1\right]^6 \to \left[0,
3\sqrt{2}\right]
\quad {\rm by} \quad 
f(\bx) = \sum_{i = 1}^3 \sqrt{x_{i1}^2 + x_{i2}^2} \ .$$
For this functional $f$,
we run our Markov chain
5000 separate times,
each from the fixed initial state
$$
\bx_0 \ = \ (0.5, 0.5, 0.5, 0.5, 0.5, 0.5)
\, ,
$$
for 500 iterations each.
We then estimate $\E[f(X_i)]$ by the average $\overline{f(X_i)}$ of the values
of the functional after $i$ iterations, averaged over the 5000 separate chains.
Since we have proven that the total variation distance is less than 0.01
after 163 iterations, the averages after 500 iterations are good estimates
of the stationary value, so we estimate using our simulations that
$$
\E_{\pi}[f] \ \approx \ \overline{f(X_{500})} \ \approx \ 2.23959
\, .
$$
On the other hand, after $i=30$ iterations, we estimate that
$$
\E[f(X_{30})] \ \approx \
\overline{f(X_{30})} \ \approx \ 2.27200
.
$$
Then, since the range of $f$
is $[0,3\sqrt{2}]$, we can estimate the total variation distance
(based on this one functional $f$) by
$$
\frac{1}{3\sqrt{2}}\left| E_{\pi}[f]  - \overline{f(X_{30})} \right|
\ \approx \
0.007638 \ < \ 0.01
.
$$
This suggests that, based on the functional $f$ at least, the
chain has approximately converged after 30 iterations.
Figure~2 shows the estimated total variation distance based on $f$
over different numbers of iterations.

\begin{figure}[ht]
    \centering
    \includegraphics[scale=0.4]{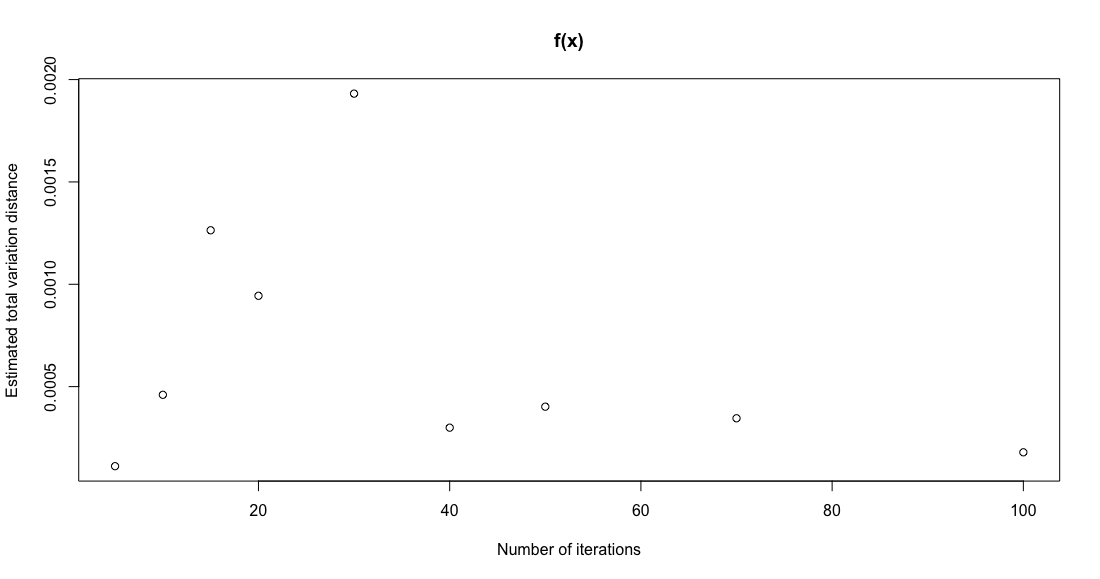}
    \caption{Estimated total variation distance based on the functional $f$,
	for the model of Sections~2 and~6.1,
	versus number of Markov chain iterations.}
\end{figure}

We also consider the following additional test functionals:
$$
g: \left[0, 1\right]^6 \to \left[0, 1\right]
\quad {\rm by} \quad
g(\bx) = x_{11}
;
$$
$$
h: \left[0, 1\right]^6 \to \left[0, \sqrt{2}\right],
\quad {\rm by} \quad
h(\bx) = ||(x_{11}, x_{12}) - (x_{21}, x_{22})||
;
$$
$$
p: \left[0, 1\right]^6
\to \left[1, e^{\sqrt{2}}\right]
\quad {\rm by} \quad
p(\bx) = \exp( ||(x_{31}, x_{32})||)
;
$$
$$
\ell: \left[0, 1\right]^6 \to \left[0, \sqrt{2}\right]
\quad {\rm by} \quad
\ell(\bx)
= \max \big(
||(x_{11}, x_{12})||, ||(x_{21}, x_{22})|| ,||(x_{31}, x_{32})|| \big)
.
$$
The estimated total variation distances based on
each of these four functionals,
as a function of the number of Markov chain iterations,
are displayed in Figure~3.
These results again
suggest that total variation distance is already below 0.01 after
just 30 iterations (though they do not show this conclusively since
the total variation distances requires a supremum over {\it all}
functionals).

\begin{figure}[ht]
    \centering
    \includegraphics[scale=0.4]{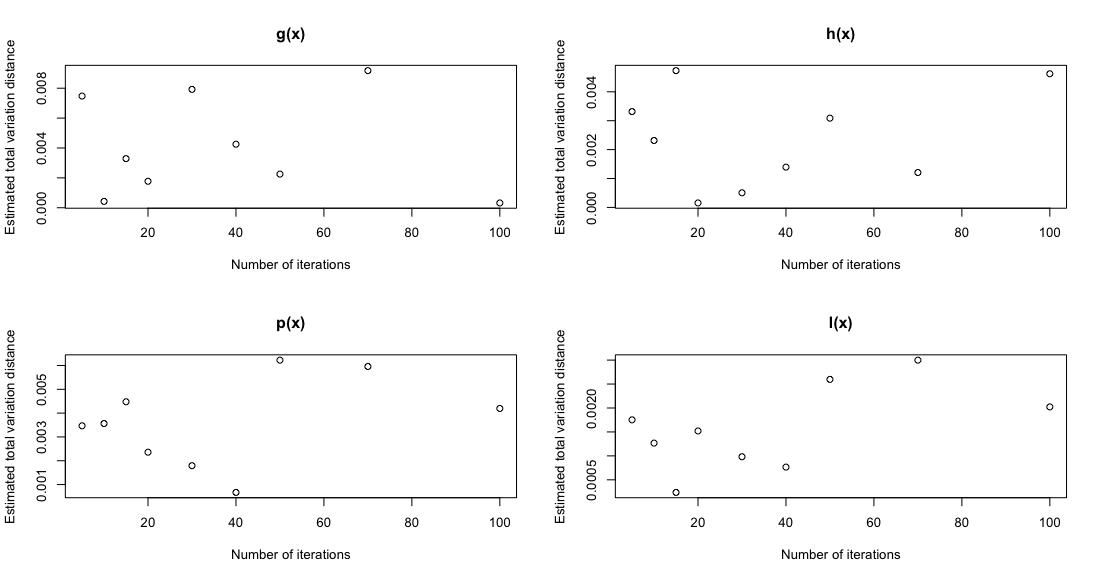}
    \caption{Estimated total variation distance based on the functionals
	$g$, $h$, $p$, and $\ell$,
	for the model of Sections~2 and~6.1,
	versus the number of Markov chain iterations.}
 \end{figure}

\subsection{One particle in $\mathbb{R}^2$}

We now consider the model of Section~3, i.e.\
a Metropolis-Hastings algorithm 
with proposals uniform on $B_x$, for
the unnormalised density~\eqref{onepi} on $\IR^2$ with one particle.
We again apply different functionals to estimate the total variation
distance.  We first let $f: \mathbb{R}^2 \to [0, 1]$ by
$f(x) = \exp(-||x||)$.  We compute by numerical integration that
$$
\E_{\pi}[f]
\ = \ { \int_{\mathbb{R}^2} f(x) \, \pi(x) \, dx
\over \int_{\mathbb{R}^2} \pi(x) \, dx }
\ \approx \ \frac{0.486}{3.189}
\ = \ 0.15240
.
$$
Similar to the previous section, we run 3000 separate chains each with initial
state $x_0 = (1, 0)$, each for 300 iterations. We then compute the mean of
$f(X_{300})$ over the 3000 chains,
and use it to estimate the total variation distance after 300 iterations
to be:
$$ \left| E_{\pi}[f]  - \overline{f(X_{300})} \right|
\ \approx \ \big| 0.15240- 0.14978\big|
\ = \ 0.00262 < 0.01
.
$$
This suggests that, based on the functional $f$ at least, the
chain has approximately converged after 300 iterations.
Figure~4 shows the estimated total variation distance based on $f$
over different numbers of iterations.

\begin{figure}[ht]
    \centering
    \includegraphics[scale=0.4]{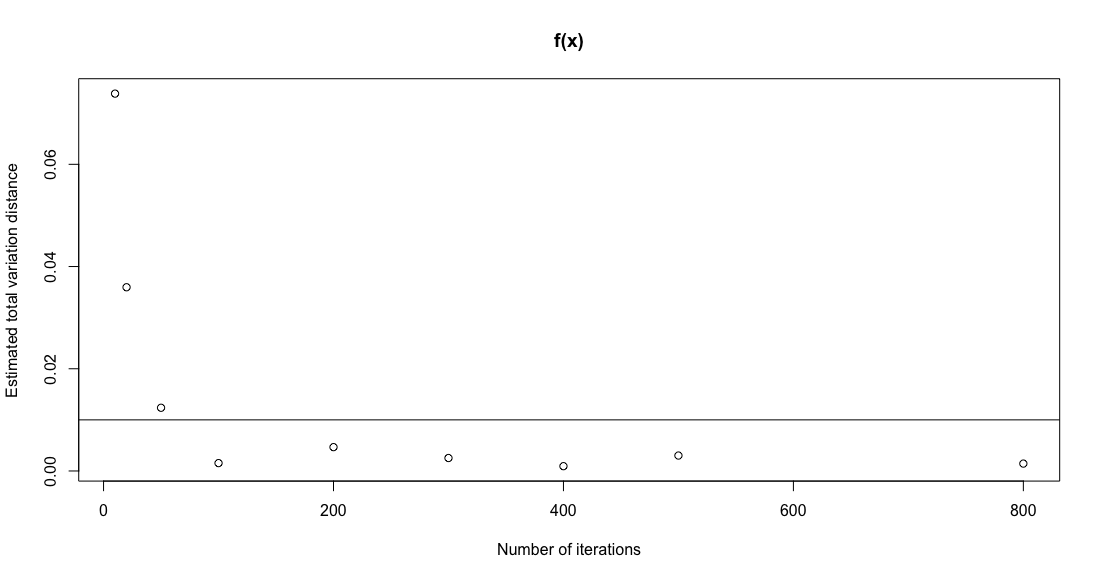}
    \caption{Estimated total variation distance based on the functional $f$,
	for the model of Sections~3 and~6.2,
	versus number of Markov chain iterations.}
\end{figure}

As before, we also consider some other test functionals.
Let
$$
g: \mathbb{R}^2 \to [0, 1]
\quad {\rm by} \quad
g(x) = \frac{x_1^2}{||x||^2}
;
$$
$$
h: \mathbb{R}^2 \to [0, 1]
\quad {\rm by} \quad
h(x) = \min\{1, \frac{1}{||x||}\}
;
$$
$$
p: \mathbb{R}^2 \to [0, 1]
\quad {\rm by} \quad
p(x) = \min\{1, |x_1|\}
;
$$
$$
\ell: \mathbb{R}^2 \to [-1, 1]
\quad {\rm by} \quad
\ell(x) = \sin(||x||)
.
$$
As before, we can use each of these functionals to estimate the total
variation distance to stationarity after different numbers of iterations,
as shown in Figure~5.  The plots suggest that total variation
distance according to each of these functionals
is below 0.01 after 300 iterations.

\begin{figure}[ht]
    \centering
    \includegraphics[scale=0.4]{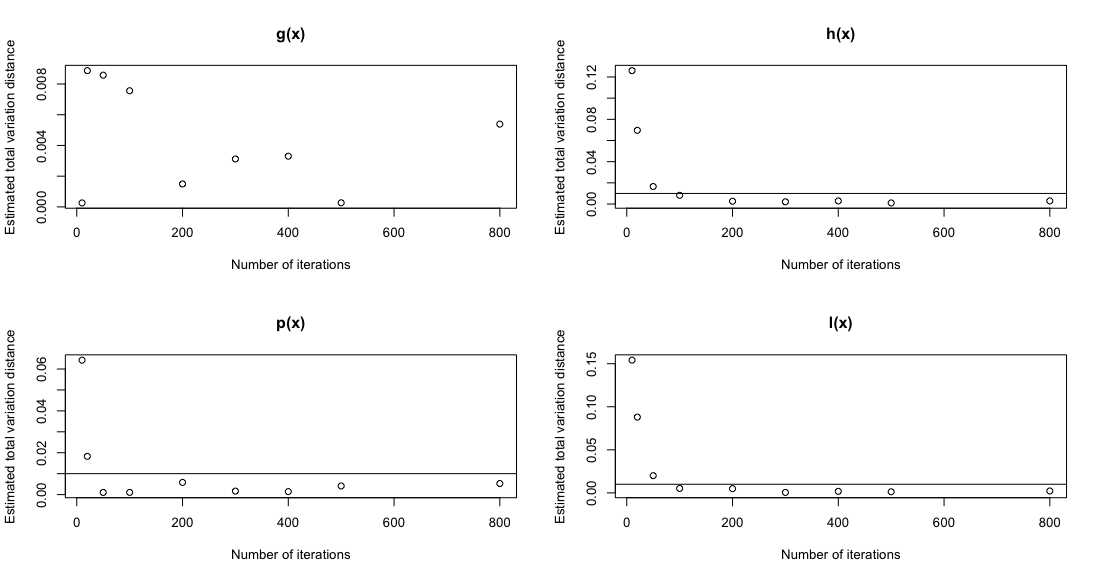}
    \caption{Estimated total variation distance based on the functionals
	$g$, $h$, $p$, and $\ell$,
	for the model of Sections~3 and~6.2,
	versus number of Markov chain iterations.}
\end{figure}

\bigskip

In summary,
both MCMC convergence
diagnostic tools and total variation distance estimation suggest
that the chains of Section~2 and Section~3 both converge
significantly more quickly than the theoretical upper
bounds derived in Sections~2 and~4.  This is not surprising, since
theoretical convergence bounds tend to be very conservative.
However, as discussed above, there is benefit in having guaranteed
theoretical convergence bounds rather
than just suggestive computer simulations which might not
accurately measure the chain's true convergence.

\section{Theorem Proofs}
\label{sec-proofs}

In this section, we prove all of the previously-stated results.

\subsection{Proof of Theorem~1}

Let
$$
\mathcal{X}' = \{(x_1, x_2, x_3) \in \mathcal{X}: \forall 1 \le i < j \le 3,
\|x_i - x_j\| \geq 1/4\}.
$$
(The value ``1/4'' is used so that $\X'$ still includes most of the mass of $\X$, but avoids states where two particles are very close thus making $\|x_i-x_j\|^{-1}$ extremely large.)
Since $\mathcal{X}'$ is compact, and $\pi(\cdot)$ is continuous and positive on $\mathcal{X'}$, therefore $\pi$
must achieve its minimum ratio $m:= \min_{\bx, \by \in \mathcal{X'}}\frac{\pi(\by)}{\pi(\bx)} > 0$ on $\X'$.
Then for any $\bx = (x_1, x_2, x_3) \in [0,1]^6$ and
measurable $A \subseteq \X$,
\begin{align*}
    P(\bx, A) &= \int_A P(\bx, d\by)\\
    &\geq \int_A P_1((x_1, x_2, x_3), dy_1) \, P_2((y_1, x_2, x_3), dy_2) \, P_3((y_1, y_2, x_3), dy_3) \\
    &\geq \int_{A \cap \mathcal{X}'} \min\bigg[1, \frac{\pi(y_1, x_2, x_3)}{\pi(x_1, x_2, x_3)}\bigg]  \min\bigg[1, \frac{\pi(y_1, y_2, x_3)}{\pi(y_1, x_2, x_3)}\bigg]  \min\bigg[1, \frac{\pi(y_1, y_2, y_3)}{\pi(y_1, y_2, x_3)}\bigg] dy,
\end{align*}
  where $P_1((x_1, x_2, x_3), B) = P((x_1, x_2, x_3), B \times \{x_2\} \times \{x_3\})$ for any measurable $B \subset [0,1]^2$ (and similarly for $P_2$
and $P_3$). Denote the three acceptance probabilities by $\alpha_1$, $\alpha_2$, $\alpha_3$ respectively.

If $\alpha_1 = 1$, then $\alpha_1\alpha_2\alpha_3 = \alpha_2\alpha_3 \geq m^2$. 
Similarly, if $\alpha_2 = 1$ or $\alpha_3=1$, then again
$\alpha_1\alpha_2\alpha_3 \geq m^2$. 
On the other hand, if  $\alpha_i < 1$ for $i = 1, 2, 3$, then $\alpha_1\alpha_2\alpha_3 = \frac{\pi(y_1, y_2, y_3)}{\pi(x_1, x_2, x_3)}\geq m \geq m^2$ (since $m \le 1$). So $$P(x, A) \ \geq \ \int_{A \cap \mathcal{X}'} m^2 \, dx \ = \ m^2 \, \Leb(A \cap \mathcal{X}'),$$ 
where $\Leb$ is Lebesgue measure on $\mathbb{R}^2$.
It follows that our algorithm satisfies a uniform minorization condition, with
$\epsilon = m^2\, \Leb(\mathcal{X}')$ and $Q(A) = \frac{\Leb(A \cap \mathcal{X}')}{\Leb(\mathcal{X}')}$.  Hence, by Proposition~1, this chain is uniformly ergodic. 

To obtain a quantitative bound, we need to compute $m^2$ and $\Leb(\mathcal{X}')$. For any $\bx \in \mathcal{X'}$,
we must have $0 \le |x_i| \le \sqrt{2}$
and $1/4 \le |x_i - x_j| \le \sqrt{2}$, thus
\begin{align*}
    0 &\le \sum_i |x_i| \le 3\sqrt{2} \text{, and } \frac{3}{\sqrt{2}} \le \sum_{i<j}|x_i - x_j|^{-1} \le 12 .
\end{align*}
Then
$$
m = \frac{\min_{\X'} \pi(\cdot)}{\max_{\X'} \pi(\cdot)} \ge \frac{e^{-c_1(3
\sqrt{2}) -c_2(12)}}{e^{-c_1(0) -c_2(3/\sqrt{2})}} = e^{-c_1(3\sqrt{2}) -
c_2(12 - 3/\sqrt{2})}.$$ Thus $$m^2 \geq \left(e^{-c_1(3\sqrt{2}) - c_2(12 -
3/\sqrt{2})}\right)^2 \geq e^{-c_1(8.49) - c_2(19.76)} .
$$
Lastly we need to compute $\Leb(\mathcal{X}')$. To make $(x_1, x_2,
x_3) \in \mathcal{X}'$, we can choose any $x_1 \in [0,1]^2$ (with area 1),
then any $x_2 \in [0, 1]^2 \setminus B(x_2, 1/4)$ (with area $\geq 1-
3.14(1/4)^2$), then any  $x_3$  $[0, 1]^2 \setminus (B(x_1, 1/4) \cup
B(x_2, 1/4))$(with area $\geq 1- 3.14(1/4)^2 - 3.14(1/4)^2$). Hence
$$\Leb(\mathcal{X}') \geq (1)\left(1 - \frac{\pi}{16}\right)
		\left(1 - \frac{\pi}{8}\right) \geq 0.48.
$$
Therefore
$$
\epsilon \ = \ m^2\Leb(\mathcal{X'})
\ \geq \ (0.48)e^{-c_1(8.49) - c_2(19.76)}.
\eqno\qed
$$

\subsection{Proof of Theorem 2 (a)}

Recall that
$$
\alpha(x, y)
\ = \ \min\bigg\{1, \frac{e^{H(x)}r_x}{e^{H(y)}r_y}\bigg\}
\ = \ \min\bigg\{1, \frac{f(r_x)}{f(r_y)}\bigg\},
$$
where $f: \mathbb{R} \to \mathbb{R}$ by $f(x) = xe^{x + \frac{1}{x}}$.
We then have
$$
f'(x) = e^{x + \frac{1}{x}} + x(1 - \frac{1}{x^2})e^{x + \frac{1}{x}} = (x - \frac{1}{x} + 1)e^{x + \frac{1}{x}},$$
so that
$$f'(x) = 0 \ \Longleftrightarrow \ x = \frac{-1 \pm \sqrt{5}}{2},
$$
and hence
$f(x)$ is decreasing on $(0, \frac{\sqrt{5} - 1}{2})$ and increasing on $(\frac{\sqrt{5} - 1}{2}, \infty)$.

Next, let
$$C_1 = \{x \in C: 1/4 \le r_x \le 2\},$$
$$C_2 =  \{x \in C: 2 \le r_x \le 4\},$$
$$D = \{x \in \mathbb{R}^2, 2 \le r_x \le 9/4\},$$
$$E_1 = \{x \in \mathbb{R}^2, 1 \le r_x \le 5/4\},$$
and $$E_2 = \{x \in \mathbb{R}^2, 3 \le r_x \le 13/4\}.$$
We shall show that $P^2(x, \cdot)$ has an overlap on $D$ for all $x \in C$. In particular, we will consider the case when the $x$ first jumps into $E_1$ and
then enters $D$ for $x \in C_1$ (similarly for $E_2$).

We know $$\alpha(x, y) =\min\left\{1, \frac{f(r_x)}{f(r_y)} \right\},$$and we have shown $f$ takes its minimum at $\frac{\sqrt{5} - 1}{2}$ and is increasing on ($\frac{\sqrt{5} - 1}{2}, \infty)$. Therefore $$m_1 := \min_{C_1 \times E_1 }\alpha(x, y) = \frac{f(\frac{\sqrt{5} - 1}{2})}{f(\frac{5}{4})} \geq 0.59, \ m_2 := \min_{C_2 \times E_2 }\alpha(x, y) = \frac{f(2)}{f(\frac{13}{4})} \geq 0.21, $$ $$  m_1' :=\min_{E_1 \times D} \alpha(x, y) = \frac{f(1)}{f(\frac{9}{4})} \geq 0.22, \ m_2':=\min_{E_2 \times D} \alpha(x, y) = \min\{ \frac{f(3)}{f(\frac{9}{4})}, 1\} =1. $$ 

For any $x \in C_1$, $y \in D$, take $M_y = \{z \in \mathbb{R}^2, r_y - 1 \le r_z \le 5/4\} \subset E_1$. Then for any $z \in M_y$, $r_z \le 5/4 \le r_x + 1$ and $r_z \geq 2 - 1 =1 \geq |r_x - 1|$. Thus $M_y \subset B_x$, and then
\begin{equation*}\begin{split}P^2(x, dy) &= \int_{B_x}P(x, dz)P(z, dy) \geq \int_{M_y}P(x, dz)P(z, dy) \\&= \frac{1}{4\pi|x|}\bigg(\int_{M_y}\alpha(x, z)\alpha(z, y)q(y, z) dz\bigg)dy  \\&\geq  \frac{1}{8\pi}\bigg(\int_{M_y} m_1m_1'  \cdot \frac{1}{4\pi|z|} dz\bigg)dy \\&= \frac{m_1m_1'}{8\pi}\bigg(2\pi \int_{r_y - 1}^{\frac{5}{4}} \frac{1}{4\pi}dr \bigg)dy  \\&= \frac{m_2m_2'}{16\pi}(\frac{9}{4} - r_y) dy \geq \frac{0.13}{16 \pi}(\frac{9}{4} - r_y)dy.\end{split}\end{equation*}

For any $x \in C_2$, $y \in D$, take $N_y = \{z \in \mathbb{R}^2, 3 \le r_z \le r_y + 1\} \subset E_2$. Similarly we have
\begin{equation*}\begin{split}P^2(x, dy) &= \int_{B_x}P(x, dz)P(z, dy) \geq \int_{M_y}P(x, dz)P(z, dy) \\&= \frac{1}{4\pi|x|}\bigg(\int_{M_y}\alpha(x, z)\alpha(z, y)q(y, z) dz\bigg)dy  \\&\geq  \frac{1}{16\pi}\bigg(\int_{M_y} m_2m_2'  \cdot \frac{1}{4\pi|z|} dz\bigg)dy \\ & =\frac{m_1m_1'}{32\pi}(r_y - 2) dy  \geq \frac{0.1}{16 \pi}( r_y - 2)dy.\end{split}\end{equation*}
Then $$P^2(x, dy) \geq 1_D \frac{1}{16\pi} \min\left\{0.13(\frac{9}{4} - ||y||), 0.1(||y||-2)\right\}dy,$$ where the size $\epsilon \geq 3.5 * 10^{-5}$.
\qed

\subsection{Proof of Theorem 2 (b)}

Since $H$ and $V$ only depend on $r_x$, we will regard them as
functions of $r_x \in \mathbb{R}$ throughout this proof.
We consider three different cases.

\medskip\noindent \textbf{Case 1}: $r_x > 4$. 

Then $r_y > 4 - 1 > \frac{\sqrt{5} - 1}{2}$ for any $y \in B_x$. So $f$ is increasing on $(r_x - 1, r_x + 1)$. For any $y \in B_x$, we have $\alpha(x, y) = 1 $ if and only if $r_y \le r_x$. Let $A_x = B(0, r_x) \setminus B(0, r_x - 1)$ (the inner part of the annulus).  Then \begin{equation*}\begin{split} PV(x) &= \int_{R^2} V(y)P(x, dy) \\&= \frac{1}{4\pi r_x}(\int_{A_x} V(y)dy + \int_{B_x \setminus A_x}V(y) \frac{f(r_x)}{f(r_y)} + \int_{B_x \setminus A_x}V(x)(1- \frac{f(r_x)}{f(r_y)}) dy \\&= \frac{1}{4\pi r_x}( \int_{A_x} V(y)dy   +    \int_{B_x \setminus A_x} (V(x) + (V(y) - V(x))\frac{f(r_x)}{f(r_y)})dy   ).\end{split}\end{equation*}
Let \begin{equation*}\begin{split}I(x, y) &=  V(x) + (V(y) - V(x))\frac{f(r_x)}{f(r_y)} = V(x) (1 + (\frac{V(y)}{V(x)} - 1)\frac{f(r_x)}{f(r_y)}) \\&= V(x)(1 + (e^{\frac{1}{2}(H(y) - H(x))} - 1)\frac{e^{H(x)}r_x}{e^{H(y)}r_y} ). \end{split}\end{equation*}
Let $u = H(y) - H(x)$, and set
$$I(x, y) = V(x)(1 + (e^{\frac{1}{2} u } - 1)e^{-u} \frac{r_x}{r_y}).$$
Then $$\int_{B_x \setminus A_x} I(x, y)dy = V(x) \left( \int_{B_x \setminus A_x}dy + \int_{B_x \setminus A_x} (e^{\frac{1}{2} u } - 1)e^{-u} \frac{r_x}{r_y}dy \right)$$ $$ = V(x)(vol(B_x \setminus A_x) + r_x \int_{B_x \setminus A_x} (e^{-\frac{1}{2} u } - e^{-u}) \frac{1}{r_y}dy ).$$
Since $u$ is a function of $r_y$ (i.e. $u$ only depends on the magnitude of $y$), $$\int_{B_x \setminus A_x} (e^{-\frac{1}{2} u } - e^{-u}) \frac{1}{r_y}dy = \int_{0}^{2\pi}\int_{r_x}^{r_x + 1}  (e^{-\frac{1}{2} u } - e^{-u}) \frac{1}{r} r dr d\theta = 2\pi \int_{r_x}^{r_x + 1} (e^{-\frac{1}{2} u } - e^{-u})  dr.$$ Since $r_x \le r_y \le r_x + 1$, $u = H(y) - H(x) = r_y - r_x + \frac{1}{r_y} - \frac{1}{r_x} \le r_y - r_x$ $\le 1$.
Note that $(e^{-\frac{1}{2} u } - e^{-u}) $ is increasing for $u \in (0, 1)$.
So $$ \int_{r_x}^{r_x + 1} (e^{-\frac{1}{2} u } - e^{-u})  dr \le \int_{r_x}^{r_x + 1}(e^{-\frac{1}{2} (r-r_x) } - e^{-(r - r_x)})dr$$
$$= \int_0^1 (e^{-\frac{1}{2}t} - e^{-t})dt = 1 + e^{-1} - 2e^{-\frac{1}{2}}. $$
Denote $(1 + e^{-1} - 2e^{-\frac{1}{2}}) $ by $m_1$. Then $$  \int_{B_x \setminus A_x} I(x, y)dy \le V(x)(vol(B_x \setminus A_x) + 2 \pi m_1 r_x) = 2\pi V(x)(r_x + \frac{1}{2} + m_1r_x).$$
(since $vol(B_x \setminus A_x) = \pi(r_x + 1)^2 - \pi r_x^2 = \pi(2r_x + 1))$.
Now consider the other part.  $$\int_{A_x} V(y)dy = 2\pi \int_{r_x - 1}^{r_x}e^{\frac{1}{2}(r + \frac{1}{r})}r dr = 2\pi V(x) \int_{r_x - 1}^{r_x}e^{\frac{1}{2}(r - r_x + \frac{1}{r} -\frac{1}{r_x})}r dr.$$
Note $$r - r_x + \frac{1}{r} -\frac{1}{r_x} = r - r_x  + \frac{r_x - r}{r r_x} \le r - r_x  + \frac{r_x - r}{12}  = \frac{11}{12}(r - r_x).$$ (the inequality follows from the fact that $rr_x \geq (r_x - 1)r_x \geq (4 - 1)4 = 12$).
So $$   \int_{A_x} V(y)dy  \le 2 \pi V(x) \int_{r_x - 1}^{r_x} e^{\frac{11}{24}(r - r_x)}rdr = 2\pi V(x) \int_{-1}^0 e^{\frac{11}{24}t}(t + r_x)dt$$   $$=2\pi V(x)  (\int_{-1}^0 te^{\frac{11}{24}t} dt + r_x\int_{-1}^0 e^{\frac{11}{24}t}dt) = 2\pi V(x)(\frac{840e^{-\frac{11}{24}} - 576}{121} + \frac{24(1 - e^{-\frac{11}{24}})}{11}r_x).$$
Denote this by $2\pi V(x)(m_2r_x + m_3)$. Then \begin{equation*} \begin{split}PV(x) &\le \frac{1}{4\pi r_x} (   2\pi V(x)(r_x + \frac{1}{2} + m_1r_x)   +    2\pi V(x)(m_2r_x + m_3)   ) \\&= \frac{V(x)}{2}(1 + \frac{1}{2r_x} + m_1 + m_2  + \frac{m_3}{r_x} )\\&\le \frac{1}{2}(1 + \frac{1}{8} + m_1 + m_2 + \frac{m_3}{4} )V(x) \text{ (as $r_x > 4$)}\\ &< 0.995V(x). \end{split}\end{equation*}

\medskip\noindent \textbf{Case 2:} $r_x < 1/4$.

In this case $|r_x - 1| = 1 - r_x > 1 - 1/4 = 3/4$, and $(r_x + 1) < 1/4 + 1 = 5/4$. So $B_x \subset (B(0, \frac{5}{4}) \setminus B(0, \frac{3}{4}) )$. Note $$\max_{y \in B_x}H(y) \le \max\{H(\frac{3}{4}), H(\frac{5}{4})\} = \max\{\frac{3}{4} + \frac{4}{3}, \frac{4}{5} + \frac{5}{4}\} = \frac{25}{12}.$$And $$H(x) \geq \frac{1}{4} + 4 = \frac{17}{4}.$$
So for any $y \in B_x$, $$V(y)/ V(x) = e^{\frac{1}{2}(H(y) - H(x))} \le e^{\frac{1}{2}(\frac{25}{12} - \frac{17}{4})} = e^{-\frac{13}{12}}.$$

Then we will show the acceptance rate is always 1. Recall $$\alpha(x, y) = \min\{1, \frac{\pi_u(y)q(y, x)}{\pi_u(x)q(x, y)} \} = \min\{1, \frac{f(r_x)}{f(r_y)}\}.$$
We showed $f(x)$ is decreasing on $(0, \frac{\sqrt{5} - 1}{2})$ and is increasing on $(\frac{\sqrt{5} - 1}{2}, \infty)$. Since $\frac{1}{4} < \frac{\sqrt{5} - 1}{2} < \frac{3}{4} $, we have $$f(r_x) \geq f(\frac{1}{4}) = \frac{1}{4}e^{\frac{17}{4}}, \  \ f(r_y) \le f(\frac{5}{4}) = \frac{5}{4}e^{\frac{41}{20}}. $$ So $$\frac{f(r_x)}{f(r_y)} \geq \frac{ \frac{1}{4}e^{\frac{17}{4}}}{ \frac{5}{4}e^{\frac{41}{20}} } > \frac{e^2}{5}> 1, \ y \in B_x. $$ Therefore $$PV(x) = \int_{B_x}q(x, y) V(y)dy \le \int_{B_x}q(x, y)e^{-\frac{13}{12}} V(x)dy = e^{-\frac{13}{12}} V(x) < 0.995V(x).$$

\medskip\noindent \textbf{Case 3:} $r_x \in [1/4, 4]$ (i.e., $x \in C$).

Let $E = B(0, \frac{1}{4})$. Note $r_y \le 5$ for all $y \in B_x$. If $y \notin E$ is proposed, since $1/4 \le r_y \le 5$ and $V(1/4) = V(4) \le V(5)$, $$V(X_{n + 1}) \le \max \{V(x), V(5)\} \le \max\{V(4), V(5)\} = e^{\frac{13}{5}}.$$

If $y \in E$ is proposed, first note this requires $|r_x - 1| < \frac{1}{4}$. So $r_x \in [\frac{3}{4}, \frac{5}{4}]$. Then $$f(r_x) \le f(\frac{5}{4}) = \frac{5}{4}e^{\frac{41}{20}}, \ f(r_y) \geq f(\frac{1}{4}) = \frac{1}{4} e^{\frac{17}{4}}.$$ So $$\frac{f(r_x)}{f(r_y)} \le \frac{\frac{5}{4}e^{\frac{41}{20}} }{  \frac{1}{4} e^{\frac{17}{4}} } < 1, \ y \in E.$$
This implies $$\alpha(x, y) = \frac{f(r_x)}{f(r_y)} < 1,\  y \in E \cap B_x.$$
Note \begin{equation*}\begin{split}PV(x) &= \int_{B_x \cap E} V(y)P(x, dy) + \int_{B_x \setminus E} V(y)P(x, dy) . \end{split}\end{equation*}
Clearly if $r_x \notin [\frac{3}{4}, \frac{5}{4}]$ (i.e. $B_x \cap E =\emptyset$ ), then $PV(x) \le V(5) = e^{\frac{13}{5}} = e^{2.6} < e^{2.7}$. Otherwise \begin{equation*}\begin{split}PV(x) &\le \int_{B_x \cap E} q(x, y) \frac{f(r_x)}{f(r_y)} V(y)dy + \int_{B_x \setminus E} V(5)P(x, dy)   \\&= \frac{f(r_x)}{4\pi r_x}\int_{B_x \cap E} \frac{e^{\frac{1}{2}H(y)}}{e^{H(y)}r_y}dy + V(5) \\&\le  \frac{f(r_x)}{4\pi r_x} 2\pi \int_{0}^\frac{1}{4} e^{-\frac{1}{2} (r + \frac{1}{r})}dr + V(5)\\&<\frac{e^{(r_x + \frac{1}{r_x)}}0.001}{2} + V(5) \\&\le \frac{e^{4 + \frac{1}{4}}}{2000} + V(5) < e^{2.7}.\end{split} \end{equation*}
Therefore
$$
PV(x) \le e^{2.7}, \quad x\in C.
$$
On the other hand, we always have $V(x) = e^{\half H(x)}
\ge e^{\half(1 + {1 \over 1})} = e$.  So, for $x\in C$,
$$
PV(x) \le e^{2.7} \le 0.995 V(x) + (e^{2.7}-0.995) \mathbf{1}_C.
\eqno\qed
$$

\subsection{Proof of Theorem~3}

It was shown in \cite{asmussen1992,aldous1993} that if two copies
$\{X_k\}_{k = 0}^{\infty}$ and $\{X_k'\}_{k = 0}^{\infty}$ of a
time-inhomogeneous Markov chain have shift-coupling times $T$ and $T'$,
then the total variation distance between the ergodic averages of their
distributions can be bounded as:
\begin{equation}\label{shiftbound}
    \left\| \frac{1}{n} \sum_{k=1}^n \P(X_k \in \cdot)
	- \frac{1}{n} \sum_{k=1}^n \P(X'_k \in \cdot) \right\|
\ \leq \ \frac{1}{n} \, \E\Big[\min(\max(T,T'), \, n)\Big].
\end{equation}
Thus, Theorem~3 will follow by constructing
copies $\{X_k\}$ and $\{X_k'\}$, with the latter in stationarity, in such
a way that we can bound these shift-coupling tail probabilities.
To do this, we generalize the construction of $\{X_k\}$ and $\{X_k'\}$ from
Section~3 of \cite{shift} to the case $n_0>1$.

Specifically, we proceed as follows.  We begin by choosing $X_0 \sim \nu$
and $X'_0 \sim \pi$ independently, and also
generate an independent random variable $W \sim Q(\cdot)$.
Then, whenever $V(X_{n})\le d$, we flip
an independent coin with probability of heads equal to $\epsilon$. If the coin
comes up heads, we set $X_{n+n_0} = W$ and $T=n+n_0$.
If the coin comes up tails, we instead generate
$X_{n+n_0} \sim \frac{1}{1-\epsilon}
(P(X_{n},\cdot)-\epsilon \, Q(\cdot))$, i.e.\ from the {\it residual}
distribution.
For completeness we then also ``fill in'' the values $X_{n + 1}, \dots, X_{n +
n_0 - 1}$ by conditional probability, according to the Markov chain
transition probabilities conditional on the already-constructed values of
$X_n$ and $X_{n+n_0}$.
If instead $V(X_n) > d$, then we simply choose $X_{n+1} \sim P(X_n,\cdot)$
as usual.
We continue this way until time $T$, i.e.\ until we get heads and set $X_T=W$.

We construct $\{X'_n\}$ and $T'$ similarly, by flipping an independent
$\epsilon$-coin whenever
$V(X'_{n})\le d$, and setting either $X'_{n + n_0} = W$ or $X'_{n+n_0} \sim
\frac{1}{1-\epsilon} (P(X'_{n},\cdot)-\epsilon \, Q(\cdot))$
(and again we ``fill in''
$X'_{n + 1}, \dots, X'_{n + n_0 - 1}$ by conditional probability),
up until the first head upon which we set $X'_{n+n_0} = W$ and $T'=n+n_0$.

This construction guarantees that $X_T = X'_{T'} = W \sim Q(\cdot)$.
We then continue the two chains identically from $W$ onwards, by
choosing $X_{T+n} = X'_{T'+n} \sim P(X_{T+n-1},\cdot)$
for $n=1,2,3,\ldots$.  Our construction ensures that each of
$\{X_n\}$ and $\{X'_n\}$ each marginally follow the transition
probabilities $P(\cdot,\cdot)$, and also that $X_{T+n} = X'_{T'+n}$ for
$n=0,1,2,\ldots$.

Now, combining the inequality~\eqref{shiftbound} with
the assumption that $P(X'_k \in \cdot) = \pi(\cdot)$ and
the standard fact
(see e.g.\ Proposition~A.2.1 of \cite{spbook}) that $\E(Z)
= \sum_{k=1}^\infty \P(Z \ge k)$
for any non-negative integer-valued random variable $Z$,
and noting that
$\P\left(\min[\max(T, T'), \, n] \geq k\right) \le
\P\left(\max(T, T') \geq k\right)$,
yields the bound:
\begin{equation}\label{sumbound}
    \left\| \frac{1}{n} \sum_{k=1}^n \P(X_k \in \cdot)
	- \pi(\cdot) \right\|
\ \leq \ \frac{1}{n} \, \sum_{k = 1}^\infty
			\P\left(\max(T, T') \geq k\right).
\end{equation}
We now bound $\P\left(\max(T, T') \geq k\right)$
for any non-negative integer $k$. 
Let $t_1,t_2,\ldots$ be the times at which we flipped a coin for
$\{X_n\}$, i.e.\ the times when $V(X_n) \le d$
excluding the ``fill in'' times.
Then, let $$N_k = \max\{i: t_i \le k\}$$
be the number coin flips up to time $k$.
Since each coin-flip yields probability $\epsilon$ of
reaching $T$ after $n_0$ additional steps, we have
for any integer $j \ge 1$ that
$\P(T \geq k, N_{k-n_0} \ge j) \le (1 - \epsilon)^j$.  Hence,
$$
\P(T \geq k)
\ = \ \P(T \geq k, \, N_{k-n_0} \geq j ) + \P(T \geq k, \, N_{k-n_0} < j)
$$
\begin{equation}\label{Tbound}
\ \le \ (1 - \epsilon)^j + \P(N_{k-n_0} < j).
\end{equation}
Then since $\lambda<1$, we have by Markov's inequality that
$$
\P(N_{k-n_0}<j)
\ = \ \P\left(t_j>{k-n_0} \right)
\ = \ \P\left(\lambda^{-t_j}>\lambda^{-k-n_0}\right)
\ \le \ \lambda^{k-n_0} \, \E\left[\lambda^{-t_j}\right].
$$
To continue, let $\tau_1 = t_1$ and $\tau_i = t_{i} - t_{i - 1}$
for $i \geq 2$.  Then by Lemma~1 below,
$$
\lambda^{k-n_0} \, \E\left[\lambda^{-t_j}\right]
\ = \ \lambda^{k-n_0} \, \E\left[\lambda^{-(\tau_1+...\tau_j)}\right]
\ \le \ \lambda^k \, \E[V(X_0)] \ (\lambda^{-n_0} A)^{j-1}.
$$
Hence, from~\eqref{Tbound},
$$
P(T\ge k) \ \le \ (1-\epsilon)^{[j]}+\lambda^{k-n_0(j-1)}A^{j-1}\E_\nu(V).
$$
Similarly,
$$
P(T'\ge k)\le (1-\epsilon)^{[j]}+\lambda^{k-n_0(j-1)} A^{j-1} \E_\pi(V).
$$
By Lemma~2 below, we have $\E_{\pi}(V) \le \frac{b}{1-\lambda}$.
Hence,
$$
\P(\max(T, T')\ge k) \ \le \ \P(T\ge k) + \P(T'\ge k)
\qquad\qquad
$$
$$
\qquad\qquad
\le
2(1-\epsilon)^{[j]}+\lambda^{k-n_0(j-1)}A^{j-1}
	\bigg(\E_\nu(V)+\frac{b}{1-\lambda}\bigg).
$$

Finally, choosing $j = \floor{rk+1} \ge rk$ and using \eqref{sumbound},
$$
    \left\| \frac{1}{n} \sum_{k=1}^n \P(X_k \in \cdot)
	- \pi(\cdot) \right\|
\qquad\qquad\qquad\qquad\qquad\qquad
$$
$$
\qquad\qquad
\ \leq \ \frac{1}{n} \, \sum_{k = 1}^\infty
	\left[ 2(1-\epsilon)^{rk}
	+ \lambda^{-n_0} (\lambda^{1-n_0 r} A^r )^k
        \bigg(\E_\nu(V)+\frac{b}{1-\lambda}\bigg) \right].
$$
Since $(1-\epsilon)^r < 1$ and $\lambda^{(1-n_0r)} A^r < 1$,
the right hand side is a geometric sum which is equal to
the claimed bound.
\qed
\bigskip

The above proof requires two lemmas.  The first is a bound on expected
values using a non-increasing expectation property, i.e.\ a partial
supermartingale argument (similar to Lemma~4 of \cite{minor1995}):

\paragraph{Lemma 1.} In the above proof of Theorem~3,
\hfil\break\bf (a) \rm\
$\E\left[\lambda^{-\tau_1}\right]
\le \E\left[V(X_0)\right]$, and
\hfil\break\bf (b) \rm\
for $i \geq 2$,
$\E\left[\lambda^{-\tau_i} | \tau_1, \dots, \tau_{i - 1}\right]
\le \lambda^{-n_0} A$.

\begin{proof} Let
$$g_i(k) = \begin{cases}\lambda^{-k} \, V(X_k), & k \le t_i \\ 0, & k >
t_i\end{cases}$$

For (a), we know that $X_k \notin C$
for any $k < t_1$, so the drift
condition implies that $g_1(k)$ has non-increasing expectation
as a function of $k$, and hence
$$\E\left[\lambda^{-\tau_1}\right] =
\E\left[\lambda^{-t_1}\right]  \le \E\left[\lambda^{-t_1}V(X_{t_1})\right]
= \E\left[g_1(t_1)\right]
\le \E\left[g_1(0)\right]  =  \E\left[V(X_0)\right].$$

For (b), for any $i \geq 2$ we know that
$X_k \notin C$ if $t_{i - 1} + n_0 \le k < t_i$, so
the drift condition implies that $g_i(k)$ has non-increasing
expectation as a function of k for $k \geq t_{i - 1} + n_0$. Hence,
\begin{equation*}\begin{split}
\E\left[\lambda^{-\tau_i} | X_{t_{i -
1}}\right]  &= \E\left[\lambda^{-(t_i - t_{i-1})} | X_{t_{i - 1}}\right]
\\& \le \E\left[\lambda^{t_{i-1}}\lambda^{-t_{i}}V(X_{t_i})  | X_{t_{i
- 1}}\right]   \\&=  \E\left[\lambda^{ t_{i-1}}g_{i}(t_i) | X_{t_{i -
1}}\right] \\ &\le \E\left[\lambda^{ t_{i-1}}g_{i}(t_{i- 1} + n_0)
| X_{t_{i - 1} }\right]  \\&\le \lambda^{ -n_0}\E\left[V(t_{i- 1} +
n_0) | X_{t_{i - 1} }\right] \\& \le \lambda^{ -n_0} \sup_{x \in C}
\E\left[V(X_1) | X_0 = x\right].
\qedhere
\end{split} \end{equation*}
\end{proof}

\bigskip
We also require a lemma which bounds $\pi(V)$, i.e.\ $\E_\pi(V)$.

\paragraph{Lemma~2:}
Suppose a $\phi$-irreducible Markov chain
on a state space $\mathcal{X}$, with
transition probabilities $P(\cdot, \cdot)$
and stationary distribution $\pi(\cdot)$, satisfies
the drift condition~\eqref{driftcond} for some function $V$
and subset $C$ and constants $\lambda<1$ and $b<\infty$.
Then the expected value of $V$ with respect to the
distribution $\pi$ satisfies the inequality $\E_\pi(V) \le b/(1-\lambda)$.

\begin{proof}
The drift condition~\eqref{driftcond} implies that our chain
satisfies \cite[Theorem~14.0.1, condition~(iii)]{MeynTweedie},
with the choice $f(x) = (1-\lambda) \, V(x)$.
Then, \cite[Theorem~14.0.1, condition~(i)]{MeynTweedie}
implies that $\E_\pi(f)<\infty$, i.e.\ $\E_\pi[(1-\lambda) \, V] < \infty$.
(The result \cite[Theorem~14.0.1]{MeynTweedie}
is actually stated assuming aperiodicity,
but it still holds in the periodic case by passing to
the lazy chain $\widebar{P} = \half(I+P)$, which is $\phi$-irreducible
and aperiodic, and has the same stationary distribution $\pi(\cdot)$,
and satisfies the drift condition~\eqref{driftcond}
with the constants $\widebar{b}=b/2$ and $\widebar{\lambda}=(1+\lambda)/2$.)
It follows that $\E_\pi(V)<\infty$.

On the other hand,
\eqref{driftcond} implies that $PV \le \lambda V + b$.
Since $\E_\pi(V)<\infty$ and $\pi \, P = \pi$,
we can take expected values with respect
to $\pi$ of both sides of this inequality to conclude that
$\E_\pi(V) \le \lambda \, \E_\pi(V) + b$.
Hence, $(1-\lambda) \, \E_\pi(V) \le b$,
so $\E_\pi(V) \le b/(1-\lambda)$, as claimed.
\end{proof}

\paragraph{Remark:} Lemma~2 can also be derived from Theorem~14.3.7 of
\cite{MeynTweedie}, with the choices
$f(x) = (1-\lambda) \, V(x)$, and $s(x) = b$, after verifying that the
chain is positive recurrent using their Theorem~14.0.1.

\subsection{Proof of Theorem~4}

For any measurable subset $S$,
\begin{align*}
|F_n(S) - \pi(S)| 
   &=\Big|\E[\text{fraction of time from 1 to n that the chain is in S}]- \pi(S)\Big|\\
   &=\bigg|\E\Big[\frac{1}{n}\sum_{k=1}^n \mathbf{1}_{X_k\in S}\Big]-\pi(S)\bigg|
   =  \bigg|\frac{1}{n}\sum_{k=1}^n \P({X_k\in S})-\pi(S)\bigg|.
\end{align*}
Thus
$$
\sup_{S\subseteq\X} |F_n(S) - \pi(S)| = \sup_S \Big|\frac{1}{n}\sum_{k=1}^n \P(X_k\in S)-\pi(S)\Big| = \Big\|\frac{1}{n}\sum_{k=1}^n \P(X_k\in \cdot)-\pi(\cdot)\Big\|,
$$
by definition of total variation distance. Also, by the triangle inequality,
$$
\bigg\|\frac{1}{n}\sum_{k=1}^n \P(X_k\in \cdot)-\pi(\cdot)\bigg\|   \le \frac{1}{n}\sum_{k=1}^n \big\|\P(X_k \in \cdot) - \pi(\cdot)\big\|.
$$
This completes the proof.
\qed

\bigskip\bigskip\noindent\bf Acknowledgements. \rm
We thank the editor and referee for very helpful suggestions which have
led to many improvements of this paper.

\bibliography{attract}{}
\bibliographystyle{plain}

\end{document}